\documentclass[11pt]{amsart}
\usepackage[parfill]{parskip}
\usepackage{tikz}
\usetikzlibrary{calc}
\usepackage{amssymb, amsmath, graphicx}
\usepackage[colorlinks=true,citecolor=black,linkcolor=black]{hyperref}

\newtheorem{theorem}{Theorem}

\newtheorem{lemma}{Lemma}
\theoremstyle{definition}
\newtheorem{definition}{Definition}
\newtheorem*{remark}{Remark}

\definecolor{pnk}{RGB}{245,169,184}
\definecolor{blu}{RGB}{91,206,250}

\title{Ordered trees and the Geode}
\author{Fern Gossow}
\date{\today}

\begin{document}

\begin{abstract}
In recent work of Wildberger and Rubine, it is shown that the formal power series $\mathbf{S}$ in the variables $t_1,t_2,\dots$ satisfying $\mathbf{S}=1+\sum_{n\geq 1} t_n\mathbf{S}^n$ has a factorisation $\mathbf{S}=1+(t_1+t_2+\cdots)\mathbf{G}$, where $\mathbf{G}$ is a power series with nonnegative coefficients called the Geode. In this note we give a combinatorial interpretation for the coefficients of $\mathbf{G}$ based on ordered trees. This amends the statement of a disproved conjecture of Wildberger and Rubine which suggests a similar (but incorrect) interpretation.
\end{abstract}

\maketitle

\section{Introduction}

For a sequence of nonnegative integers $\mathbf{m}=(m_1,m_2,\dots)$ with finite sum, define the \emph{hyper-Catalan number}
\[C_\mathbf{m}=\frac{(m_1+2m_2+3m_3+\cdots)!}{(1+m_2+2m_3+\cdots)!\,m_1!\,m_2!\,m_3!\cdots}.\]
For variables $t_1,t_2,\dots$ set $\mathbf{t}^\mathbf{m}:=t_1^{m_1}t_2^{m_2}\cdots$ and define $\mathbf{S}:=\sum_{\mathbf{m}}C_\mathbf{m} \mathbf{t}^\mathbf{m}$, which is the unique formal power series satisfying the functional equation
\[\mathbf{S}=1+\sum_{n\geq 1}t_n\mathbf{S}^n.\]
This series has garnered recent interest due to Wildberger and Rubine \cite{wr}, who use it to find roots of arbitrary integer polynomials. Their paper also contains a brief history of the hyper-Catalan numbers.

\begin{remark}
The definitions of $C_\mathbf{m}$ and $\mathbf{S}$ in \cite{wr} assume that $m_1=0$, but our generalised setup can be easily obtained \cite{g}.
\end{remark}

Combinatorial interpretations of $C_\mathbf{m}$ date back to at least 1940 in work of Etherington \cite{ee}. We overview these constructions in Section \ref{sect:S} by defining families of sets $\{S_\mathbf{m}\}$ and $\{T_\mathbf{m}\}$ based on polygonal dissections and ordered trees respectively such that
\[C_\mathbf{m}=\vert S_\mathbf{m}\vert=\vert T_\mathbf{m}\vert\]
for every $\mathbf{m}$. In the case $\mathbf{m}=(0,n,0,\dots)$, $C_\mathbf{m}$ is the $n^\text{th}$ Catalan number, $S_\mathbf{m}$ is the set of triangulations of the $(n+2)$-gon, and $T_\mathbf{m}$ is the set of binary trees with $n+1$ leaves.

In their paper, Wildberger and Rubine observe that $\mathbf{S}$ has a surprising factorisation
\[\mathbf{S}=1+(t_2+t_3+\cdots)\mathbf{G}\]
where $\mathbf{G}$ is another power series with nonnegative coefficients. They name $\mathbf{G}$ \emph{the Geode} and give a number of conjectures about its coefficients. Many of these were soon proven \cite{az}\cite{r1}, but a general formula is still unknown. Recent work of Gessel \cite{g} shows that the coefficients of $\mathbf{G}$ count certain lattice paths.

It was conjectured in \cite{wr} that the coefficient of $t_2^{m_2}t_3^{m_3}$ in $\mathbf{G}$ counts ordered trees (trees with a designated root node and an ordering on the children of any node) with $m_2$ nodes with two children, $m_3$ nodes with three children, and a single additional leaf node. However, this conjecture is false, and refined attempts in follow-up work have been so far unsuccessful \cite{r2}. The main objective of this note is to provide a version of the construction which correctly interprets the coefficients of $\mathbf{G}$.

\begin{theorem}\label{thm:G}
For a sequence $\mathbf{m}=(m_1,m_2,\dots)$ of nonnegative integers with finite sum, let $L_\mathbf{m}$ denote the number of nodes $v$ among ordered trees $T$ such that:
\begin{itemize}
\item $T$ has $m_n$ nodes with exactly $n$ children for every $n\geq 1$,
\item $v$ is a leaf node of $T$, and
\item $v$ is visited before any non-leaf node in post-order traversal.
\end{itemize}
Then $\mathbf{G}=\sum_\mathbf{m} L_\mathbf{m}\mathbf{t}^\mathbf{m}$.
\end{theorem}

For the following ordered trees $T$, the nodes $v$ which satisfy the above conditions have been circled.

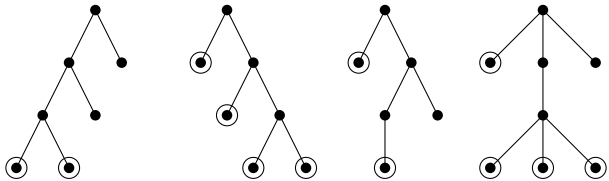
\[\begin{tikzpicture}[scale = 0.35]
\fill [black] (0,8) circle (0.2);
\fill [black] (-1,6) circle (0.2);
\fill [black] (1,6) circle (0.2);
\fill [black] (-2,4) circle (0.2);
\fill [black] (0,4) circle (0.2);
\fill [black] (-3,2) circle (0.2);
\fill [black] (-1,2) circle (0.2);
\draw (0,8) -- (-3,2);
\draw (-2,4) -- (-1,2);
\draw (-1,6) -- (0,4);
\draw (0,8) -- (1,6);
\draw [black] (-3,2) circle (0.4);
\draw [black] (-1,2) circle (0.4);

\fill [black] (5,8) circle (0.2);
\fill [black] (4,6) circle (0.2);
\fill [black] (6,6) circle (0.2);
\fill [black] (5,4) circle (0.2);
\fill [black] (7,4) circle (0.2);
\fill [black] (6,2) circle (0.2);
\fill [black] (8,2) circle (0.2);
\draw (5,8) -- (8,2);
\draw (7,4) -- (6,2);
\draw (6,6) -- (5,4);
\draw (5,8) -- (4,6);
\draw [black] (6,2) circle (0.4);
\draw [black] (8,2) circle (0.4);
\draw [black] (5,4) circle (0.4);
\draw [black] (4,6) circle (0.4);

\fill [black] (11,8) circle (0.2);
\fill [black] (10,6) circle (0.2);
\fill [black] (12,6) circle (0.2);
\fill [black] (11,4) circle (0.2);
\fill [black] (13,4) circle (0.2);
\fill [black] (11,2) circle (0.2);
\draw (10,6) -- (11,8) -- (12,6) -- (11,4) -- (11,2);
\draw (12,6) -- (13,4);
\draw [black] (10,6) circle (0.4);
\draw [black] (11,2) circle (0.4);

\fill [black] (17,8) circle (0.2);
\fill [black] (15,6) circle (0.2);
\fill [black] (17,6) circle (0.2);
\fill [black] (19,6) circle (0.2);
\fill [black] (17,4) circle (0.2);
\fill [black] (15,2) circle (0.2);
\fill [black] (17,2) circle (0.2);
\fill [black] (19,2) circle (0.2);
\draw (15,6) -- (17,8) -- (19,6);
\draw (15,2) -- (17,4) -- (19,2);
\draw (17,8) -- (17,2);
\draw [black] (15,6) circle (0.4);
\draw [black] (15,2) circle (0.4);
\draw [black] (17,2) circle (0.4);
\draw [black] (19,2) circle (0.4);

\end{tikzpicture}\]

\section{What does \texorpdfstring{$\mathbf{S}$}{S} count?}\label{sect:S}

In this section we define the hyper-Catalan series $\mathbf{S}$ and give two combinatorial interpretations for its coefficients. These results are essentially due to Etherington \cite{ee} and can be found explicitly in more recent work of Schuetz and Gwyn \cite{s}, but having a self-contained proof will be useful when we construct novel interpretations for the coefficients of $\mathbf{G}$ in the next section.

Let $\mathbf{m}=(m_1,m_2,\dots)$ denote a sequence of nonnegative integers with finite sum. We set the notation $\mathbf{0}:=(0,0,\dots)$ and $\mathbf{e}_n:=(\dots,0,1,0,\dots)$ with a $1$ in the $n^\text{th}$ place. Let $t_1,t_2,\dots$ be variables and define $\mathbf{t}^\mathbf{m}:=t_1^{m_1}t_2^{m_2}\cdots$. Let $\mathbf{S}$ be the unique formal power series in the variables $t_1,t_2,\dots$ such that
\[\mathbf{S}=1+\sum_{n\geq 1}t_n\mathbf{S}^n.\]

\begin{definition}
A \emph{subdigon} is a polygon with a designated edge called the \emph{roof}, which has been subdivided into faces given by noncrossing arcs between its vertices. The \emph{type} of a subdigon is the sequence $\mathbf{m}=(m_1,m_2,\dots)$, where $m_1$ counts the number of bigon faces, $m_2$ the number of triangle faces, $m_3$ the number of quadrilateral faces and so on.
\end{definition}

\begin{remark}
As previously noted, this is a slight generalisation of the definition in \cite{wr}, since it allows for multiple edges between the same vertices. Setting $m_1=0$ recovers the original definition, but we don't want to let bigons be bygones.
\end{remark}

We draw subdigons with their roof dashed. Below is a subdigon of type $\mathbf{m}=(2,3,2,1,0,\dots)$.

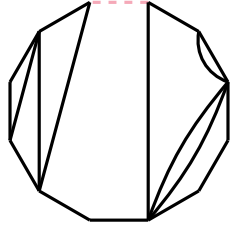
\[\begin{tikzpicture}[rotate = 75, scale = 0.3]
\foreach \i in {1,...,11} {
	\draw [black, very thick] (\i*30:5) -- (\i*30+30:5);
}
\draw [pnk, dashed, very thick] (0:5) -- (30:5);

\draw [black, very thick] (1*30:5) -- (5*30:5);
\draw [black, very thick] (2*30:5) -- (5*30:5);
\draw [black, very thick] (2*30:5) -- (4*30:5);
\draw [black, very thick] (0*30:5) -- (7*30:5);
\draw [black, very thick] (7*30:5) to [bend left = 10] (10*30:5);
\draw [black, very thick] (7*30:5) to [bend right = 10] (10*30:5);
\draw [black, very thick] (10*30:5) to [bend left = 40] (11*30:5);
\end{tikzpicture}\]

\begin{definition}
An \emph{ordered tree} is a tree graph with a designated root node, where the children of a node are assigned an ordering. An ordered tree has \emph{type} $\mathbf{m}=(m_1,m_2,\dots)$ if $m_1$ nodes have one child, $m_2$ nodes have two children and so on. This is also known as the \emph{downdegree sequence} \cite{s}.
\end{definition}

We draw ordered trees with their root at the top, and the children of any node below their parent, arranged left-to-right according to their ordering. Below is an ordered tree of type $\mathbf{m}=(2,3,2,1,0,\dots)$.

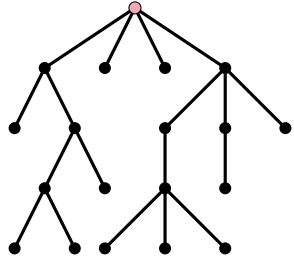
\[\begin{tikzpicture}[scale = 0.4]
\draw [black, very thick] (0,4) -- (-3,2);
\draw [black, very thick] (0,4) -- (-1,2);
\draw [black, very thick] (0,4) -- (1,2);
\draw [black, very thick] (0,4) -- (3,2);

\draw [black, very thick] (-3,2) -- (-4,0);
\draw [black, very thick] (-3,2) -- (-2,0);
\draw [black, very thick] (3,2) -- (1,0);
\draw [black, very thick] (3,2) -- (3,0);
\draw [black, very thick] (3,2) -- (5,0);

\draw [black, very thick] (-2,0) -- (-3,-2);
\draw [black, very thick] (-2,0) -- (-1,-2);
\draw [black, very thick] (1,0) -- (1,-2);
\draw [black, very thick] (3,0) -- (3,-2);

\draw [black, very thick] (-3,-2) -- (-4,-4);
\draw [black, very thick] (-3,-2) -- (-2,-4);
\draw [black, very thick] (1,-2) -- (-1,-4);
\draw [black, very thick] (1,-2) -- (1,-4);
\draw [black, very thick] (1,-2) -- (3,-4);

\draw [black, fill = pnk] (0,4) circle (0.2);
\fill [black] (-3,2) circle (0.2);
\fill [black] (-1,2) circle (0.2);
\fill [black] (1,2) circle (0.2);
\fill [black] (3,2) circle (0.2);

\fill [black] (-4,0) circle (0.2);
\fill [black] (-2,0) circle (0.2);
\fill [black] (1,0) circle (0.2);
\fill [black] (3,0) circle (0.2);
\fill [black] (5,0) circle (0.2);

\fill [black] (-3,-2) circle (0.2);
\fill [black] (-1,-2) circle (0.2);
\fill [black] (1,-2) circle (0.2);
\fill [black] (3,-2) circle (0.2);

\fill [black] (-4,-4) circle (0.2);
\fill [black] (-2,-4) circle (0.2);
\fill [black] (-1,-4) circle (0.2);
\fill [black] (1,-4) circle (0.2);
\fill [black] (3,-4) circle (0.2);
\end{tikzpicture}\]

Let $S_\mathbf{m}$ and $T_\mathbf{m}$ denote the set of subdigons and ordered trees of type $\mathbf{m}$ respectively. It is known that these sets have the same size, and this size matches the coefficient of $\mathbf{t}^\mathbf{m}$ in $\mathbf{S}$. In the proof below we construct a bijection between subdigons and ordered trees, which will match together the above examples. (Can you guess the map?)

\begin{theorem}\label{thm:S} There is an equality of formal power series
\[\sum_{\mathbf{m}}\vert S_\mathbf{m}\vert t^\mathbf{m}=\sum_{\mathbf{m}}\vert T_\mathbf{m}\vert t^\mathbf{m}=\mathbf{S}.\]
\end{theorem}

\begin{proof}
\emph{First equality}. We construct an explicit bijection $S_\mathbf{m}\to T_\mathbf{m}$ for every $\mathbf{m}$. Given a subdigon, define the \emph{central face} as the one containing the roof. Create a tree consisting of a root node (corresponding to the roof) with children $v_1,\dots,v_n$ corresponding to the other edges $e_1,\dots,e_n$ of the central face, listed in counterclockwise order starting from the roof.

Each $v_i$ becomes the root of a subtree, which is recursively constructed using the ``subsubdigon'' who faces are those which are separated from the central face by $e_i$, and whose roof is now $e_i$. For our running example, the initial decomposition will give four subsubdigons, drawn below with their new roofs.

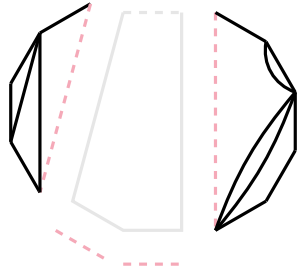
\[\begin{tikzpicture}[rotate = 75, scale = 0.3]
\draw [black!10!white, very thick] (1*30:5) -- (5*30:5) -- (6*30:5) -- (7*30:5) -- (0*30:5);

\draw [black!10!white, dashed, very thick] (1*30:5) -- (0*30:5);

\foreach \i in {1,...,4} {
	\draw [black, very thick] ($(90:1.5)+(\i*30:5)$) -- ($(90:1.5)+(\i*30+30:5)$);
}

\draw [pnk, dashed, very thick] ($(90:1.5)+(1*30:5)$) -- ($(90:1.5)+(5*30:5)$);
\draw [black, very thick] ($(90:1.5)+(2*30:5)$) -- ($(90:1.5)+(5*30:5)$);
\draw [black, very thick] ($(90:1.5)+(2*30:5)$) -- ($(90:1.5)+(4*30:5)$);

\draw [pnk, dashed, very thick] ($(5.5*30:1.5)+(5*30:5)$) -- ($(5.5*30:1.5)+(6*30:5)$);

\draw [pnk, dashed, very thick] ($(6.5*30:1.5)+(6*30:5)$) -- ($(6.5*30:1.5)+(7*30:5)$);

\draw [pnk, dashed, very thick] ($(9.5*30:1.5)+(0*30:5)$) -- ($(9.5*30:1.5)+(7*30:5)$);
\foreach \i in {7,...,11} {
	\draw [black, very thick] ($(9.5*30:1.5)+(\i*30:5)$) -- ($(9.5*30:1.5)+(\i*30+30:5)$);
}
\draw [black, very thick] ($(9.5*30:1.5)+(7*30:5)$) to [bend left = 10] ($(9.5*30:1.5)+(10*30:5)$);
\draw [black, very thick] ($(9.5*30:1.5)+(7*30:5)$) to [bend right = 10] ($(9.5*30:1.5)+(10*30:5)$);
\draw [black, very thick] ($(9.5*30:1.5)+(10*30:5)$) to [bend left = 40] ($(9.5*30:1.5)+(11*30:5)$);
\end{tikzpicture}\]

The inverse map is obtained by recursively constructing a 
subdigon for each subtree of the root's $n$ children, and gluing these to a new roofed polygon of size $n+1$ in counterclockwise order (starting from the roof, which is not glued onto). 

Under the bijection, edges of a subdigon correspond to nodes of the associated tree. Internal nodes correspond to faces, and faces with $n+1$ edges (a roof and $n$ other edges) correspond to nodes with $n$ children. Hence, subdigons of type $\mathbf{m}$ become ordered trees of type $\mathbf{m}$ as required.

\emph{Second equality}. We prove that $\mathbf{T}:=\sum_\mathbf{m}\vert T_\mathbf{m}\vert t^\mathbf{m}$ satisfies the same functional equation as $\mathbf{S}$. After observing that $\vert T_{\mathbf{0}}\vert=1$, the statement $\mathbf{T}=1+\sum_{n\geq 1}t_n\mathbf{T}^n$ is equivalent to the existence of a bijection
\[T_\mathbf{m} \longrightarrow \bigcup_{n\geq 1}\bigcup_{\genfrac{}{}{0pt}{1}{\mathbf{m}_1,\dots,\mathbf{m}_n;}{\mathbf{e}_n+\mathbf{m}_1+\cdots+\mathbf{m}_n=\mathbf{m}}} \,T_{\mathbf{m}_1}\times\cdots\times T_{\mathbf{m}_n}\]
for all $\mathbf{m}\neq\mathbf{0}$.

This bijection is given by mapping an ordered tree $\tau\in T_\mathbf{m}$ to the tuple $(\tau_1,\dots,\tau_n)$, where $(v_1,\dots,v_n)$ are the children of the root in $\tau$ (listed in order), and $\tau_i$ is the subtree rooted at $v_i$. The degree of $\tau$ is $\mathbf{e}_n$ (for the root node) plus the degrees of $\tau_1,\dots,\tau_n$, as required. The inverse map is given by joining the roots of a list of ordered trees to a new root node.
\end{proof}

\section{What does \texorpdfstring{$\mathbf{G}$}{G} count?}\label{sect:G}

In this section we prove Theorem \ref{thm:G} by constructing sets $\overline{T}_\mathbf{m}$ whose sizes are the coefficients of the formal power series $\mathbf{G}$ satisfying
\[\mathbf{S}=1+(t_1+t_2+\cdots)\mathbf{G}.\]
By Theorem \ref{thm:S} and the above equation, we expect that any $\tau\in T_\mathbf{m}$ of type $\mathbf{m}\neq \mathbf{0}$ can be canonically decomposed as a pair $(n,\overline{\tau})$ where $n\geq 1$ is an integer and $\overline{\tau}\in\overline{T}_{\mathbf{m}-\mathbf{e}_n}$. This will correspond to finding the left-most vertex $v$ in $\tau$ whose children $v_1,\dots,v_n$ are all leaves, removing those leaves, and marking $v$ to create $\overline{\tau}$. The preimage of this map corresponds to unmarking $v$ and attaching to it any positive number of children, all of which are leaves.

We also make a related subdigon-based construction, which may be useful in the search to find more combinatorial interpretations for the coefficients of $\mathbf{G}$.

\begin{definition}
A face in a subdigon is called \emph{external} if it has a single internal edge, where the roof counts as internal.
\end{definition}

\begin{definition}
An internal node in an ordered tree is called \emph{clawed} if all of its children are leaves.
\end{definition}

Under the bijection in the proof of Theorem \ref{thm:S}, we observe that external faces correspond to clawed nodes. We also observe that every subdigon has at least one external face, and every ordered tree has at least one clawed node (except for unique object of type $\mathbf{0}$ in both cases). Below we shade the external faces and clawed nodes of our running examples.

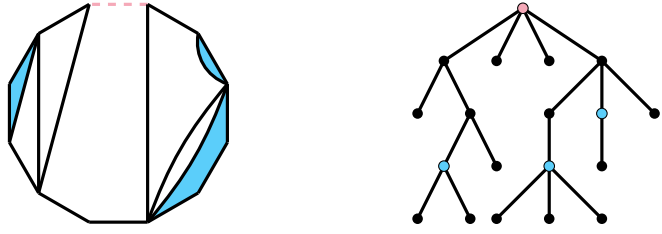
\[\vcenter{\hbox{
\begin{tikzpicture}[rotate = 75, scale = 0.3]
\fill [blu] (2*30:5) -- (3*30:5) -- (4*30:5) -- cycle;
\fill [blu] (7*30:5) to [bend right = 10] (10*30:5) -- (9*30:5) -- (8*30:5) -- cycle;
\fill [blu] (10*30:5) to [bend left = 40] (11*30:5) -- cycle;

\foreach \i in {1,...,11} {
	\draw [black, very thick] (\i*30:5) -- (\i*30+30:5);
}
\draw [pnk, dashed, very thick] (0:5) -- (30:5);

\draw [black, very thick] (1*30:5) -- (5*30:5);
\draw [black, very thick] (2*30:5) -- (5*30:5);
\draw [black, very thick] (2*30:5) -- (4*30:5);
\draw [black, very thick] (0*30:5) -- (7*30:5);
\draw [black, very thick] (7*30:5) to [bend left = 10] (10*30:5);
\draw [black, very thick] (7*30:5) to [bend right = 10] (10*30:5);
\draw [black, very thick] (10*30:5) to [bend left = 40] (11*30:5);
\end{tikzpicture}}}\qquad\qquad\qquad \vcenter{\hbox{
\begin{tikzpicture}[scale = 0.35]
\draw [black, very thick] (0,4) -- (-3,2);
\draw [black, very thick] (0,4) -- (-1,2);
\draw [black, very thick] (0,4) -- (1,2);
\draw [black, very thick] (0,4) -- (3,2);

\draw [black, very thick] (-3,2) -- (-4,0);
\draw [black, very thick] (-3,2) -- (-2,0);
\draw [black, very thick] (3,2) -- (1,0);
\draw [black, very thick] (3,2) -- (3,0);
\draw [black, very thick] (3,2) -- (5,0);

\draw [black, very thick] (-2,0) -- (-3,-2);
\draw [black, very thick] (-2,0) -- (-1,-2);
\draw [black, very thick] (1,0) -- (1,-2);
\draw [black, very thick] (3,0) -- (3,-2);

\draw [black, very thick] (-3,-2) -- (-4,-4);
\draw [black, very thick] (-3,-2) -- (-2,-4);
\draw [black, very thick] (1,-2) -- (-1,-4);
\draw [black, very thick] (1,-2) -- (1,-4);
\draw [black, very thick] (1,-2) -- (3,-4);

\draw [black, fill = pnk] (0,4) circle (0.2);
\fill [black] (-3,2) circle (0.2);
\fill [black] (-1,2) circle (0.2);
\fill [black] (1,2) circle (0.2);
\fill [black] (3,2) circle (0.2);

\fill [black] (-4,0) circle (0.2);
\fill [black] (-2,0) circle (0.2);
\fill [black] (1,0) circle (0.2);
\draw [black, fill = blu] (3,0) circle (0.2);
\fill [black] (5,0) circle (0.2);

\draw [black, fill = blu] (-3,-2) circle (0.2);
\fill [black] (-1,-2) circle (0.2);
\draw [black, fill = blu] (1,-2) circle (0.2);
\fill [black] (3,-2) circle (0.2);

\fill [black] (-4,-4) circle (0.2);
\fill [black] (-2,-4) circle (0.2);
\fill [black] (-1,-4) circle (0.2);
\fill [black] (1,-4) circle (0.2);
\fill [black] (3,-4) circle (0.2);
\end{tikzpicture}}}\]

\begin{lemma}\label{lemma:G}
Let $\overline{S}_\mathbf{m}$ denote the set of subdigons of type $\mathbf{m}$ with a marked external edge, which occurs before or inside of the first external face encountered when travelling counterclockwise from the roof. Then $\mathbf{G}=\sum_{\mathbf{m}}\vert\overline{S}_\mathbf{m}\vert t^{\mathbf{m}}$.
\end{lemma}
\begin{proof}
Let $\overline{\mathbf{S}}=\sum_{\mathbf{m}}\vert\overline{S}_\mathbf{m}\vert t^{\mathbf{m}}$. We prove that $\mathbf{S}=1+(t_1+t_2+\cdots)\overline{\mathbf{S}}$ by exhibiting a bijection
\[S_{\mathbf{m}} \to \bigcup_{\genfrac{}{}{0pt}{1}{n;}{m_n\geq 1}}\,\{n\}\times \overline{S}_{\mathbf{m}-\mathbf{e}_n}\]
for all $\mathbf{m}\neq\mathbf{0}$.

Given $\sigma\in S_\mathbf{m}$, let $f$ denote the first external face encountered when travelling counterclockwise from the roof. Consider the map $\sigma\mapsto (n,\overline{\sigma})$, where $f$ has $n+1$ edges and $\overline{\sigma}\in\overline{S}_{\mathbf{m}-\mathbf{e}_n}$ is obtained by deleting $f$ and marking the unique internal (but now external) edge. The inverse map is given by attaching a polygon with $n+1$ sides to the marked edge.
\end{proof}

We now state the analogous result for trees. By the bijection in Theorem \ref{thm:S}, we know that external edges in subdigons correspond to leaves in ordered trees. It remains to find a description for leaves whose corresponding external edges occur before or within the first external face.

\begin{definition}
The \emph{post-order traversal} of an ordered tree $\tau\in T_\mathbf{m}$ is an ordering on its nodes defined recursively as follows, beginning at the root:
\begin{itemize}
\item if this node has no children, `visit' the node (add it next in the order),
\item else run the procedure on each child of this node in order.
\end{itemize}
\end{definition}

This method of visiting nodes is also known as \emph{depth-first search}, since it visits all children of a node before visiting the node itself. Let $\overline{T}_\mathbf{m}$ denote the set of ordered trees of type $\mathbf{m}$ with a marked leaf node, which is visited before any non-leaf node during post-order traversal. Theorem \ref{thm:G} is equivalent to the claim
\[\mathbf{G}=\sum_{\mathbf{m}}\vert \overline{T}_\mathbf{m}\vert \mathbf{t}^\mathbf{m}.\]

\begin{proof}[Proof of Theorem \ref{thm:G}]
Suppose $\sigma\in S_\mathbf{m}$ maps to $\tau \in T_\mathbf{m}$ under the bijection from Theorem \ref{thm:S}. The first node in $\tau$ visited during post-order traversal corresponds to the external edge in $\sigma$ counterclockwise-adjacent to the roof. Each subsequent leaf node visited corresponds to the next-adjacent external edge. The first non-leaf node visited in $\tau$ is the clawed node corresponding to the first external face in $\sigma$ reached by travelling counterclockwise from the roof. Hence, leaf nodes in $\tau$ which occur before any non-leaf node correspond to external edges of $\sigma$ which occur before or inside the first external face. This gives $\vert \overline{S}_\mathbf{m}\vert=\vert\overline{T}_\mathbf{m}\vert$ and we apply Lemma \ref{lemma:G}.
\end{proof}

\quad {\bf Acknowledgment}. Thanks to Dean Rubine for useful suggestions which improved the exposition in this work.

\medskip

\bibliography{geode}
\bibliographystyle{plain}

\end{document}